\documentclass[12pt]{article}

\usepackage{amsfonts}

\newtheorem{theorem}{Theorem}
\newtheorem{lemma}[theorem]{Lemma}
\newtheorem{corollary}[theorem]{Corollary}

\newtheorem{proposition}[theorem]{Proposition}
\newtheorem{example}[theorem]{Example}
\newtheorem{remark}[theorem]{Remark}

\def\bit{\begin{itemize}}
\def\eit{\end{itemize}}

\def\bc{\begin{center}}
\def\ec{\end{center}}

\def\bthm{\begin{theorem}}
\def\ethm{\end{theorem}}

\def\bcor{\begin{corollary}}
\def\ecor{\end{corollary}}

\def\bprop{\begin{proposition}}
\def\eprop{\end{proposition}}

\def\blem{\begin{lemma}}
\def\elem{\end{lemma}}

\def\bex{\begin{example}}
\def\eex{\end{example}}

\def\brem{\begin{remark}}
\def\erem{\end{remark}}

\def\prf{\noindent{\bf Proof. }}

\def\bdes{\begin{description}}
\def\edes{\end{description}}

\def\iti{\item[(i)]}
\def\itii{\item[(ii)]}

\def\beq{\begin{equation}}
\def\eeq{\end{equation}}

\def\ben{\begin{enumerate}}
\def\een{\end{enumerate}}

\def\beqar{\begin{eqnarray}}
\def\eeqar{\end{eqnarray}}

\def\beqarr{\begin{eqnarray*}}
\def\eeqarr{\end{eqnarray*}}

\def \non{{\nonumber}}

\def\RR{{\mathbb R}}  

\def\SS{{\mathbb S}}

\def\cB{\mathcal{B}}

\def\cF{\mathcal{F}}

\def\qed{\vbox{\hrule\hbox{\vrule height 1.5 ex\kern 1 ex\vrule}\hrule}}

\def\P{{\mathsf P}} 

\def\E{{\mathsf E}} 

\def\eps{\varepsilon}

\def\part{\partial}

\begin{document}

\title{The noise of a Brownian sticky flow is black.}
\author{Yves Le Jan and Olivier Raimond}
\maketitle



We refer the reader to \cite{ljr2} for the construction of a sticky
flow on the circle. A Brownian sticky flow is a sticky flow whose one
point motion is a Brownian motion.

Let us first recall Tsirelson's definition of a noise (see
\cite{tsirelson1,tsirelson2})~: A {\em noise} consists of a
probability space $(\Omega,\cF,\P)$, a family
$(\cF_{s,t})_{s\leq t}$ of sub-$\sigma$-fields (also called a
factorization) of $\cF$ and a $L^2$-continuous one-parameter group
$(T_h)_{h\in \RR}$ of transformations of $\Omega$ preserving $\P$ such
that
\bdes \iti for all $s\leq t\leq u$,
$\cF_{s,u}=\sigma(\cF_{s,t}\cup\cF_{t,u})$.
\itii for all $t_1<\cdots <t_n$, $(\cF_{t_{i-1},t_i})_{1\leq i\leq n}$
is a family of independent $\sigma$-fields.
\itii for all $s\leq t$ and $h\in\RR$, $T_h(\cF_{s,t})=\cF_{s+h,t+h}$.
\edes
In \cite{ljr1}, a noise $N=(\Omega,\cF,(\cF_{s,t}),\P,(T_h))$ is
associated to a stochastic flow of kernels on a locally compact metric
space $M$ (more precisely to the law of a stochastic flow of kernels),
with $\Omega=\prod_{s\leq t} E$ ($E$ is the space of kernels on $E$),
$\cF=\otimes_{s<t}\cB(E)$ ($\cB(E)$ is the Borel $\sigma$-field on
$E$), $K_{s,t}(\omega)=\omega(s,t)$, $\P$ is the law of the stochastic
flow of kernels and $T_h$ is defined by
$T_h(\omega)(s,t)=\omega(s+h,t+h)$. We call this noise  
{\em the noise of the stochastic flow of kernels} $K$.

A process $X=(X_{s,t})_{s\leq t}$ is called a 
{\em centered linear representation} of a noise 
$N=(\Omega,\cF,(\cF_{s,t}),\P,(T_h))$ ($X$ is also called a decomposable
process) if for all $s\leq t$, $X_{s,t}\in L^2_0(\cF_{s,t})$ ($L^2_0$
is the set of all $L^2$-functions with mean 0) and for all $s\leq t\leq u$,
a.s. $X_{s,u}=X_{s,t}+X_{t,u}$. We denote by $H_0^{\hbox{lin}}$ the space
of centered linear representations of $N$. The {\em linear part} of
$N$, denoted $N^{\hbox{lin}}=(\Omega,\cF,(\cF^{\hbox{lin}}_{s,t}),\P,(T_h))$, is a
subnoise of $N$ (i.e. $\cF^{\hbox{lin}}_{s,t}\subset\cF_{s,t}$) where
$\cF_{s,t}^{\hbox{lin}}=\sigma(X_{u,v},~X\in H_0^{\hbox{lin}},~s\leq
u\leq v\leq t)$. A noise is called {\em black} if $N^{\hbox{lin}}$ is
a trivial noise, or equivalently if $H_0^{\hbox{lin}}=\{0\}$.

\medskip
The purpose of this note is to prove
\bthm The noise of a Brownian sticky flow is black. \ethm

A noise $N$ is called {\em continuous} (see \cite{tsirelson2}) (or the
 factorization $(\cF_{s,t})$ is called continuous) if for every $s<t$,
$\cup_{\eps>0}\cF_{s+\eps,t-\eps}$ generates $\cF_{s,t}$ and
 $\cup_{n=1}^\infty\cF_{-n,n}$ generates $\cF$.
\blem The noise of a stochastic flow of kernels is continuous. \elem
\prf It is enough to show that for all
$s\leq t$ and all $X\in L^2(\P)$, 
\beq\label{continue}
\lim_{\eps\to 0+}\E[X|\cF_{s+\eps,t-\eps}]=\E[X|\cF_{s,t}] \quad
\hbox{ in } L^2(\P).\eeq

Let us remark that $L^2(\P)$ is the closure of the vector space
spanned by functions of the form
\beq \label{form}
X=\prod_{i=1}^n\left(\prod_{j=1}^{k_i} K_{t_{i-1},t_i}f^i_j(x^i_j)\right),\eeq
where $f^i_j$ are Lipschitz continuous functions on $M$,
$x_j^i\in M$ and $t_0<\cdots<t_n$. In the case $X=K_{0,1}f(x)$,
(\ref{continue})
 is satisfied for all $0<s<t<1$ since (in the following,
$\P_t$ denotes the Feller semigroup of the one-point motion)
$$\E[K_{0,1}f(x)|\cF_{s+\eps,t-\eps}] =
\P_{s+\eps}K_{s+\eps,t-\eps}\P_{1-(t-\eps)}f(x)$$
which converges in $L^2(\P)$ towards
$\P_{s}K_{s,t}\P_{1-t}f(x)=\E[K_{0,1}f(x)|\cF_{s,t}]$ as $\eps\to
0+$. Similarly, it can be shown that (\ref{continue}) is satisfied for
all $X$ of the form (\ref{form}) and all $s\leq t$. Thus $N$ is
continuous. \qed

\medskip
From now on, $N$ denotes the noise of a sticky flow.

\smallskip
Let $H_1(s,t)=\{X_{s,t}:~X\in H_0^{\hbox{lin}}\}$ be the first chaos
of $N$. For every $X\in L^2_0(\Omega)$, $H_{s,t}(X)$
denotes the orthogonal projection of $X$ on $H_1(s,t)$. We set
$H=H_{0,1}$. Since $N$ is continuous, $H(X)=\lim_{n\to\infty}
\sum_{k=1}^{2^n}\E[X|\cF_{(k-1)2^{-n},k2^{-n}}]$
(see proposition 6.a.2 in \cite{tsirelson2}). Note that $N$ is
black if $H(X)=0$ for all $X\in L^2_0(\P)$, or equivalently for all
$X\in L^2_0(\cF_{0,1})$.

\blem\label{essentiel}
Let $f$ be a Lipschitz function on $(\SS^1)^d$ and $h$ be a
nonnegative continuous function on $\SS^1$.
\beq H\left(\langle K_{0,1}^{\otimes d}f,h\rangle_{L^2(m_d)}
-\langle \P^{(d)}_1f,h\rangle_{L^2(m_d)}\right)=0.\eeq \elem
\prf We set $X=\langle K_{0,1}^{\otimes d}f,h\rangle_{L^2(m_d)}
-\langle \P^{(d)}_1f,h\rangle_{L^2(m_d)}$. The result holds if
$$\sum_{k=1}^{2^n}\E[X|\cF_{(k-1)2^{-n},k2^{-n}}]$$
converges towards 0 in $L^2(\P)$ as $n\to\infty$. Since the terms
in the sum are independent, this holds if
\beq\label{sum} \lim_{n\to\infty}\sum_{k=1}^{2^n}
\E\left[\left(\E[X|\cF_{(k-1)2^{-n},k2^{-n}}]\right)^2\right] = 0. \eeq
Note that for all $0\leq s<t\leq 1$, 
$$\E[\langle K_{0,1}^{\otimes d}f,h\rangle_{L^2(m_d)}|\cF_{s,t}] =
\langle \P^{(d)}_sK_{s,t}^{\otimes d}
\P^{(d)}_{1-t}f,h\rangle_{L^2(m_d)}.$$
We set $\mu_s(dx)=\int_{x_0\in
  (\SS^1)^d}h(x_0)\P^{(d)}_s(x_0,dx)m_d(dx_0)$, $g_t=\P^{(d)}_{1-t}f$
and $\eps=t-s$. Then we have
$$\E\left[\left(\E[X|\cF_{s,t}]\right)^2\right]
=\mu_s^{\otimes 2} \left(\P^{(2d)}_\eps (g_t\otimes g_t)-\P^{(d)\otimes 2}_\eps
(g_t\otimes g_t)\right),$$
where $\P^{(2d)}_t$ is the semigroup of the $2d$ point motion of the
sticky flow and $\P^{(d)\otimes 2}_t$ is the semigroup of two independent
$d$ point motions. Note that $\mu_s$ is absolutely continuous with
respect to $m_d$ and we have $\frac{d\mu_s}{dm_d}\leq
\|h\|_\infty<\infty$ (since for $f\geq 0$, $\mu_sf=\langle
\P_s^{(d)}f,h\rangle_{L^2(m_d)} = \langle
f,\P_s^{(d)}h\rangle_{L^2(m_d)} \leq \|h\|_{\infty}m_df$). We also
have that $Lip(g_t)\leq Lip(f)$ (where $Lip(f)$ denotes the Lipschitz
constant of $f$).

It is easy to check, using the compatibility relations of the family
$(\P^{(n)}_t)$ that for all $(x,y)\in (\SS^1)^d\times (\SS^1)^d$
\beqar\left(\P^{(2d)}_\eps (g_t\otimes g_t) - 
\P^{(d)\otimes 2}_\eps (g_t\otimes g_t)\right)(x,y)\hskip-100pt&&\non\\
&=& \frac{1}{2}(\E^{(d)\otimes 2}_{(x,y)}-\E^{(2d)}_{(x,y)})
[(g_t(X_\eps)-g_t(Y_\eps))^2],\label{eq}\eeqar 
where under $\P^{(2d)}_{(x,y)}$ (resp. under $\P^{(d)\otimes 2}_{(x,y)}$)
$(X_t,Y_t)_t$ is the $2d$ point motion started at $(x,y)$ (resp. $X_t$
and $Y_t$ are independent $d$ point motions respectively started at
$x$ and $y$).

\blem \label{lemunif} Let $g$ be a Lipschitz function on
$(\SS^1)^d$. Then for all positive $t$, all $x$ and $y$ in $(\SS^1)^d$,
we have
\beq \label{equnif}  
(\E^{(d)\otimes 2}_{(x,y)}-\E^{(2d)\otimes 2}_{(x,y)})
[(g(X_t)-g(Y_t))^2] \leq (4dLip(g))\times t. \eeq \elem
\prf We have
\beqarr
(\E^{(d)\otimes 2}_{(x,y)}-\E^{(2d)\otimes 2}_{(x,y)})
[(g(X_t)-g(Y_t))^2] \hskip-155pt &&\\
&=&(\E^{(d)\otimes 2}_{(x,y)}-\E^{(2d)\otimes 2}_{(x,y)})
[((g(X_t)-g(x))+(g(y)-g(Y_t)))^2]\\
&+& (\E^{(d)\otimes 2}_{(x,y)}-\E^{(2d)\otimes 2}_{(x,y)})
[(g(x)-g(y))^2]\\
&+& 2(\E^{(d)\otimes 2}_{(x,y)}-\E^{(2d)\otimes 2}_{(x,y)})
[(g(x)-g(y))((g(X_t)-g(x))+(g(y)-g(Y_t)))].\eeqarr
It is easy to see that the second and the third terms vanish. Thus
\beqarr
(\E^{(d)\otimes 2}_{(x,y)}-\E^{(2d)\otimes 2}_{(x,y)})
\left[\left(g(X_t)-g(Y_t)\right)^2\right] \hskip-80pt &&\\
&\leq&2\E^{(d)\otimes 2}_{(x,y)}
\left[\left(g(X_t)-g(x)\right)^2+\left(g(y)-g(Y_t)\right)^2\right]\\
&\leq& 2Lip(g)\E^{(d)\otimes 2}_{(x,y)}\left[d(X_t,x)^2+d(Y_t,y)^2\right]\\
&\leq& 4dLip(g)\times t. \eeqarr
This proves the lemma. \qed

\medskip Set $\tau = \inf\{t,~\{X_t\}\cap\{Y_t\}\neq\emptyset\} =
\inf_{1\leq i,j\leq d}\tau_{i,j}$, where
$\tau_{i,j}=\inf\{t,~X^i_t=Y^j_t\}$. Note that $((X_t,Y_t),~t\leq
\tau)$ has the same law under $\P^{(2d)}_{(x,y)}$ and under
$\P^{(d)\otimes 2}_{(x,y)}$.

\blem \label{lemtau} There exists a constant $C$ depending only on $d$
such that for all positive $\eps$,
\beq \label{eqtau}\P^{(d)\otimes 2}_{m_d^{\otimes 2}}[\tau<\eps]
<C\sqrt{\eps}. \eeq\elem
\prf We have $\P^{(d)\otimes 2}_{m_d^{\otimes 2}}[\tau<\eps]
\leq \sum_{1\leq i,j\leq d} 
\P^{(d)\otimes 2}_{m_d^{\otimes 2}}[\tau_{i,j}<\eps]$. Let us remark
that for all $i,j$
$$\P^{(d)\otimes 2}_{m_d^{\otimes 2}}[\tau_{i,j}<\eps] =
\frac{1}{2\pi^2}\int_{0<x<y<2\pi}\P_{y-x}[T\leq 2\eps]dx dy,$$
with $T=\inf\{t,B_t\in\{0,2\pi\}\}$ and where, under $\P_z$, $B_t$ is a
Brownian motion starting at $z$ (we use the fact that the law of
$(Y^j_{t/2}-X^i_{t/2},~t<2\tau)$ under $\P^{(d)\otimes 2}_{(x,y)}$ is
the same as the law of $(B_t,~t<T)$ under $\P_z$). 
Since $\P_z[T\leq 2\eps]\leq \P_z[T_0\leq 2\eps]+\P_{2\pi-z}[T_0\leq 2\eps]$, where
$T_0=\inf\{t,~B_t=0\}$, we have
$$\P^{(d)\otimes 2}_{m_d^{\otimes 2}}[\tau_{i,j}<\eps] \leq
  \frac{2}{\pi}\int_0^{2\pi}\P_z[T_0<2\eps]~dz.$$
Using the reflection principle, we get $\P_z[T_0<2\eps]\leq
  2\P_0[B_{2\eps}>z]$. Thus
\beqarr
\P^{(d)\otimes 2}_{m_d^{\otimes 2}}[\tau_{i,j}<\eps] &\leq&
  \frac{4}{\pi}\int_0^{2\pi}\P_0[B_{2\eps}>z]~dz\\
&\leq& \frac{4}{\pi}\sqrt{2\eps}\E_0[|B_1|]. \eeqarr
This proves the lemma. \qed

\blem\label{lemestim} There exists a constant $C$ depending only on
$d$, $Lip(f)$ and $\|h\|_\infty$ such that for all $0<s<t<1$ and $\eps=t-s$,
\beq \label{eqestim} \mu_s^{\otimes 2}
\left(\P^{(2d)}_\eps (g_t\otimes g_t)-
\P^{(d)\otimes 2}_\eps (g_t\otimes g_t)\right) \leq C\eps^{3/2}.\eeq \elem
\prf Using first the fact that $((X_t,Y_t),~t\leq
\tau)$ has the same law under $\P^{(2d)}_{(x,y)}$ and under
$\P^{(d)\otimes 2}_{(x,y)}$, then the strong Markov property at time
$\tau$ and finally lemma \ref{lemunif}, we have
\beqarr
(\E^{(d)\otimes2}_{(x,y)}-\E^{(2d)}_{(x,y)})[(g(X_\eps)-g(Y_\eps))^2]
 \hskip-100pt&&\\
&=& (\E^{(d)\otimes2}_{(x,y)}-\E^{(2d)}_{(x,y)})
    [1_{\tau<\eps}(g(X_\eps)-g(Y_\eps))^2]\\
&=& \E^{(d)\otimes 2}_{(x,y)}
    \left[1_{\tau<\eps}(\E^{(d)\otimes2}_{(X_\tau,Y_\tau)} -
      \E^{(2d)}_{(X_\tau,Y_\tau)})
	[(g(X_{\eps-\tau})-g(Y_{\eps-\tau}))^2]\right]\\
&\leq& (4dLip(f))\times\eps\times\P^{(d)\otimes 2}_{(x,y)}[\tau<\eps].
\eeqarr
Using now (\ref{eq}) and lemma \ref{lemtau}, we get
\beqarr
\mu_s^{\otimes 2}
\left(\P^{(2d)}_\eps (g_t\otimes g_t)-
\P^{(d)\otimes 2}_\eps (g_t\otimes g_t)\right) \hskip-80pt &&\\
&\leq& (2dLip(f)) \times \eps \times \int \mu_s^{\otimes 2}(dx,dy)
\P^{(d)\otimes 2}_{(x,y)}[\tau<\eps]\\
&\leq& (2dLip(f)\|h\|_\infty^2) \times \eps \times
   \P^{(d)\otimes 2}[\tau<\eps]\\
&\leq& (2dLip(f)\|h\|_\infty^2C) \times \eps^{3/2}.\qquad  \qed\eeqarr

\medskip

Lemma \ref{lemestim} permits to finish the proof of lemma
\ref{essentiel}. Indeed, lemma \ref{lemestim} implies that
$$\E[(\E[X|\cF_{(k-1)2^{-n},k2^{-n}}])^2]\leq C\times 2^{-3n/2}.$$
This implies that
$$\sum_{k=1}^{2^n}\E[(\E[X|\cF_{(k-1)2^{-n},k2^{-n}}])^2]\leq C\times
2^{-n/2},$$
which converges towards $0$ as $n\to\infty$. Thus $H(X)=0$. \qed

\medskip
\noindent{\bf Proof of the theorem.} We are going to show that
$H(X)=0$ for $X\in L^2_0(\cF_{0,1})$. Let $V$ denote the vector space spanned
by functions of the form 
\beq\label{form2}
X=\prod_{i=1}^n \langle K_{t_{i-1},t_i}^{\otimes d}f^i,h^i\rangle_{L^2(m_{d})},\eeq
for all $n\geq 1$, $d\geq1$,  $(f_i)_i$ a family of Lipschitz functions on
$(\SS^1)^d$, $(h_i)$ a family of nonnegative continuous functions on
$(\SS^1)^d$ and $0=t_0<t_1<\cdots<t_n=1$. Let $V_0$ be the set of
functions in $V$ such that $\E[X]=0$. Then $V$ is
dense in $L^2(\cF_{0,1})$ and $V_0$ is dense in
$L_0^2(\cF_{0,1})$. Thus to prove that $H(X)=0$ for all $X\in
L^2_0(\cF_{0,1})$, it is enough to prove $H(X)=0$ for all $X\in V_0$. Let
$X$ be in the form (\ref{form2}), then for all $i$, there exists a
constant $c_i$ such that  
$$\E[X|\cF_{t_{i-1},t_i}] =  c_i\langle
K_{t_{i-1},t_i}^{\otimes d}f^i,h^i\rangle_{L^2(m_d)}.$$
Thus, if we take $X$ such that for all $i$,
$\E[\langle K_{t_{i-1},t_i}^{\otimes d}f^i,h^i\rangle_{L^2(m_d)}]=0$,
$\E[X]=0$ (note that $V_0$ is spanned by functions in this form) and
\beqarr H(X)
&=& \sum_{i=1}^n H_{t_{i-1},t_i}(X)\\
&=& \sum_{i=1}^n c_i H_{t_{i-1},t_i}
     (\langle K_{t_{i-1},t_i}^{\otimes d}f^i,h^i\rangle_{L^2(m_d)}).
\eeqarr
Note that lemma \ref{essentiel} shows that $H(X)=0$.
This proves the theorem. \qed

\brem Note that this proof extends easily to prove that the noise generated by
Arratia's coalescing flow is black. This was originally proved by
Tsirelson (cf \cite{tsirelson2}). Our proof is based on the criterium
he gave to prove that a noise is black (zero quadratic variation). \erem

\bibliographystyle{apalike}

\end{document}